\documentclass[fleqn]{mat01}
\usepackage{times,mathtimy,amssymb,latexsym,amscd}
\usepackage[mathscr]{eucal}
\begin{document}

\newtheorem{theor}{\bf Theorem}
\newtheorem{lem}{Lemma}
\newtheorem{propo}{\rm PROPOSITION}

\setcounter{page}{121} \firstpage{121}

\font\xxx=eusm10 at 7.6pt

\def\CCA{\mbox{\xxx{\char'101}}}
\def\CCB{\mbox{\xxx{\char'102}}}
\def\CCU{\mbox{\xxx{\char'125}}}
\def\CCX{\mbox{\xxx{\char'130}}}
\def\CCH{\mbox{\xxx{\char'110}}}

\def\BA{{\bf A}}
\def\BB{{\bf B}}
\def\BC{{\mathbb C}}
\def\BE{{\bf E}}
\def\BF{{\bf F}}
\def\BG{{\bf G}}
\def\BH{{\mathbb H}}
\def\BJ{{\bf J}}
\def\BL{{\bf L}}
\def\BN{{\mathbb N}}
\def\BQ{{\mathbb Q}}
\def\BR{{\mathbb R}}
\def\BS{{\bf S}}
\def\BT{{\mathbb T}}
\def\BU{{\bf U}}
\def\BV{{\bf V}}
\def\BX{{\bf X}}
\def\BZ{{\mathbb Z}}
\def\ba{{\bf a}}
\def\bb{{\bf b}}
\def\bc{{\bf c}}
\def\bd{{\bf d}}
\def\be{{\bf e}}
\def\bff{{\bf f}}
\def\bg{{\bf g}}
\def\bh{{\bf h}}
\def\bi{{\bf i}}
\def\bj{{\bf j}}
\def\bk{{\bf k}}
\def\bl{{\bf l}}
\def\bm{{\bf m}}
\def\bn{{\bf n}}
\def\bp{{\bf p}}
\def\bq{{\bf q}}
\def\br{{\bf r}}
\def\bs{{\bf s}}
\def\bt{{\bf t}}
\def\bu{{\bf u}}
\def\bv{{\bf v}}
\def\bw{{\bf w}}
\def\bx{{\bf x}}
\def\by{{\bf y}}
\def\bz{{\bf z}}
\let\PAR=\S
\let\a=\alpha
\let\b=\beta
\let\g=\gamma
\let\G=\Gamma
\let\d=\delta
\let\D=\Delta
\let\e=\epsilon
\let\z=\zeta
\let\i=\iota
\let\th=\theta
\let\TH=\Theta
\let\k=\kappa
\let\l=\lambda
\let\L=\Lambda
\let\m=\mu
\let\n=\nu
\let\r=\rho
\let\s=\sigma
\let\S=\Sigma
\let\t=\tau
\let\OM=\Omega
\let\om=\omega
\let\x=\xi
\let\X=\Xi
\let\y=\eta
\def\CA{{\mathscr A}}
\def\CB{{\mathscr B}}
\def\CC{{\mathscr C}}
\def\CD{{\mathscr D}}
\def\CE{{\mathscr E}}
\def\CF{{\mathscr F}}
\def\CG{{\mathscr G}}
\def\CH{{\mathscr H}}
\def\CI{{\mathscr I}}
\def\CJ{{\mathscr J}}
\def\CK{{\mathscr K}}
\def\CL{{\mathscr L}}
\def\CM{{\mathscr M}}
\def\CN{{\mathscr N}}
\def\CO{{\mathscr O}}
\def\CP{{\mathscr P}}
\def\CQ{{\mathscr Q}}
\def\CR{{\mathscr R}}
\def\CS{{\mathscr S}}
\def\CT{{\mathscr T}}
\def\CU{{\mathscr U}}
\def\CV{{\mathscr V}}
\def\CW{{\mathscr W}}
\def\CX{{\mathscr X}}
\def\CY{{\mathscr Y}}
\def\CZ{{\mathscr Z}}
\let\Inf=\inf
\let\Sup=\sup
\let\inf=\infty
\let\sbs=\subset
\let\sps=\supset
\def\inv{{-1}}
\let\del=\partial
\def\blb{\bigl\{}
\def\brb{\bigr\}}
\def\symdif{{\mathbin\vartriangle}}
\def\card{\mathop{\rm card}}
\def\BRN{{\bf R}^{2n}}
\let\<=\langle
\let\>=\rangle
\def\tpi{{2\pi i}}
\def\tpixx{{2\pi i\x\cdot x}}
\def\tpiyy{{2\pi i\y\cdot y}}
\def\HN{{{\bf H}_n}}
\def\half{{\textstyle{\frac12}}}
\def\qed{\hbox{\quad\vrule height6pt width2.5pt}}
\def\pd#1#2{{\frac{\del#1}{\del#2}}}
\font\twelvebf=cmbx12
\font\sans=cmss10
\font\twelverm=cmr12
\font\bfi=cmbxti10
\def\intint{\int\!\!\!\int}
\def\intintint{\int\!\!\!\int\!\!\!\int}
\def\intintintint{\int\!\!\!\int\!\!\!\int\!\!\!\int}
\def\intn{\int\!\!\cdots\!\!\int}
\def\nat{{\,\natural\,}}
\let\bar=\overline
\def\Proof{{\sl Proof:}\quad}
\def\sha{{\,\sharp\,}}
\let\T=\dagger
\def\singsupp{\,\hbox{\rm sing supp}\,}
\def\smatrix#1#2#3#4{\Bigl({#1\atop#3}{#2\atop#4}\Bigr)}
\def\Im{\mathop{\rm Im}\nolimits}
\def\Re{\mathop{\rm Re}\nolimits}
\def\fourth{{\textstyle{\frac14}}}
\def\sech{\mathop{\rm sech}\nolimits}
\def\csch{\mathop{\rm csch}\nolimits}
\def\grad{\nabla}
\def\Aut{\mathop{\rm Aut}\nolimits}
\let\sq=\square
\def\sgn{\mathop{\rm sgn}\nolimits}
\def\Res{\mathop{\rm Res}\nolimits}
\def\supp{\mathop{\rm supp}\nolimits}
\let\dsp=\displaystyle
\def\BIGbreak{\goodbreak \vskip 24pt plus 4 pt minus 4 pt}
\outer\def\proclaim #1. #2\par{\medbreak
\noindent{\bf#1.\enspace}{\sl#2\par \medbreak}}
\def\sectionhead#1{\BIGbreak\noindent{\twelvebf #1}\mark{#1}
\par\nobreak\bigskip\nobreak}
\def\tr{\mathop{\rm tr}\nolimits}
\hyphenation{iso-metry}
\hyphenation{iso-met-ric}
\def\ind{\mathop{\rm ind}\nolimits}
\def\norm{{|\mskip -1.5mu |\mskip -1.5mu |}}
\def\U{{\rm sup}}
\let\emptyset=\varnothing
\def\diam{\mathop{\rm diam}}
\def\loc{{\rm loc}}
\def\CCC{{C_c^\inf}}
\def\endhat{{\hat{\mathstrut\ }}}
\def\sinc{\mathop{\rm sinc}\nolimits}
\def\hattimes{\mathbin{\hat\times}}
\def\curl{\mathop{\rm curl}}
\def\div{\mathop{\rm div}}
\def\mapright#1{\smash{\mathop{\longrightarrow}\limits^{#1}}}
\def\mapleft#1{\smash{\mathop{\longleftarrow}\limits^{#1}}}
\let\w=\wedge

\title{The abstruse meets the applicable: Some aspects\\ of time-frequency analysis}

\markboth{G~B~Folland}{Some aspects of time-frequency analysis}

\author{G~B~FOLLAND}

\address{Department of Mathematics, University of Washington,
Seattle, WA 98195-4350, USA\\
\noindent E-mail: folland@math.washington.edu}

\volume{116}

\mon{May}

\parts{2}

\pubyear{2006}

\Date{MS received 6 February 2006}

\begin{abstract}
The area of Fourier analysis connected to signal processing theory
has undergone a rapid development in the last two decades. The
aspect of this development that has received the most publicity is
the theory of wavelets and their relatives, which involves
expansions in terms of sets of functions generated from a single
function by translations and dilations. However, there has also
been much progress in the related area known as
\emph{time-frequency analysis} or \emph{Gabor analysis}, which
involves expansions in terms of sets of functions generated from a
single function by translations and modulations. In this area
there are some questions of a concrete and practical nature whose
study reveals connections with aspects of harmonic and functional
analysis that were previously considered quite pure and perhaps
rather exotic. In this expository paper, I~give a survey of some
of these interactions between the abstruse and the applicable. It
is based on the thematic lectures which I~gave at the Ninth
Discussion Meeting on Harmonic Analysis at the Harish-Chandra
Research Institute in Allahabad in\break October 2005.
\end{abstract}


\maketitle

\section{Background and statements of theorems}

In this section we develop some basic ideas of time-frequency
analysis and lead up to the main results to be discussed in this
paper. For a more comprehensive account of this subject we
recommend the book of Gr\"ochenig \cite{gro}; see also Daubechies
\cite{daub2}.

For an $L^2$ function $f$ on $\BR^d$, the Fourier transform and
its inversion formula
\begin{equation*}
\hat{f} (\om)=\int f(t)\hbox{e}^{-\tpi \om\cdot t}\,\hbox{d}
t,\qquad f(t)=\int \hat{f}(\om) \hbox{e}^{\tpi \om\cdot
t}\,\hbox{d}\om
\end{equation*}
provide the expansion of $f$ in terms of the pure `sine waves'
$\hbox{e}^{\tpi \om\cdot t}$. (When \hbox{$d=1$} it is often
appropriate to refer to $t$ as `time' and $\om$ as `frequency.' We
shall sometimes use these terms even when dealing with general
$d$, for which the basic mathematical structure is exactly the
same.) The Fourier transform is a marvelous tool, but since it
involves the~whole function $f$ at once, it is not an efficient
way to analyze the ways in which different frequencies enter into
$f$ at different times (as is of paramount importance, for
example, in music).

One way to produce a `local' Fourier analysis of a function $f$ is
to use the \emph{windowed Fourier transform}, which is also known
as the \emph{short-time Fourier transform} and is closely related
to the \emph{cross-ambiguity function} of radar theory and the
\emph{Fourier--Wigner transform} in \cite{fol1}. The idea is
simple: one fixes an\vadjust{\pagebreak} $L^2$ function $\phi$ on
$\BR^d$ such that $\|\phi\|_2=1$ and considers the Fourier
transform not of $f$ but of $f$ multiplied by translates of
$\phi$, obtaining a function $V_\phi f$ on $\BR^d\times\BR^d$:
\begin{equation}\label{1}
V_\phi f(\om,x)=\int f(t)\bar{\phi(t-x)}\hbox{e}^{-\tpi \om\cdot
t}\,\hbox{d}t.
\end{equation}
$V_\phi f$ is called the \emph{windowed Fourier transform of $f$
with window $\phi$}. (The reasons for the complex conjugation on
$\phi$ and the normalization $\|\phi\|_2=1$ will become apparent
shortly.) If one takes $\phi$ to be a nonnegative function
supported on an interval $I$ centered at the origin (or at least
negligibly small outside $I$), then $V_\phi f(\om,x)$ measures how
much the frequency $\om$ contributes to the portion of $f$ that
lives on the interval $I+x$ centered at $x$.

To obtain the inversion formula for the operator $V_\phi$, we
observe that the map $(\phi,f)\mapsto V_\phi f$ is the restriction
to functions of the form $F(x,t)=\bar{\phi(x)}f(t)$ of the linear
map $\tilde{V}\hbox{:}\ L^2(\BR^{2d})\to L^2(\BR^{2d})$ defined by
\begin{equation*}
\tilde{V}\!F(\om,x)=\int F(t-x,t)\hbox{e}^{-\tpi\om\cdot
t}\,\hbox{d}t.
\end{equation*}
$\tilde{V}$ is the composition of the measure-preserving change of
variable $(x,t)\mapsto (t-x,t)$ with the Fourier transform in the
second variable, so it is unitary on $L^2$. It follows that the
windowed Fourier transform is an isometry from $L^2(\BR^d)$ into
$L^2(\BR^{2d})$:
\begin{equation*}
\|V_\phi f\|_2=\|\phi\|_2\|f\|_2=\|f\|_2.
\end{equation*}
Therefore, we have $V^*_\phi V_\phi=I$, and an easy calculation of
$V_\phi^*$ then yields the inversion formula
\begin{equation}\label{2}
f(t)=\iint V_\phi f(\om,x)\phi(t-x)\hbox{e}^{\tpi\om\cdot
t}\,\hbox{d}x\,\hbox{d}\om.
\end{equation}
(As with the ordinary Fourier transform, this integral is
absolutely convergent only for $f$ in a dense subspace of $L^2$
and must be interpreted by a limiting process in general.)

Let us look at this from another viewpoint. For $\om,x\in\BR^d$ we
introduce the \emph{modulation} operator $M_\om$ and the
\emph{translation} operator $T_x$ by
\begin{equation*}
M_\om f(t)=\hbox{e}^{\tpi \om\cdot t}f(t),\quad T_xf(t)=f(t-x).
\end{equation*}
Thus
\begin{align}\label{3}
M_\om T_x f(t) &=\hbox{e}^{\tpi \om\cdot t}f(t-x),\nonumber\\[.3pc]
T_xM_\om f(t)&=\hbox{e}^{\tpi\om\cdot(t-x)}f(t-x)=
\hbox{e}^{-\tpi\om\cdot x}M_\om T_xf(t),\phantom{0}\hskip -1pc
\end{align}
and \eqref{1} and \eqref{2} become
\begin{equation}\label{4}
V_\phi f(\om,x)=\langle f,M_\om T_x\phi\rangle,\quad f=\iint
V_\phi f(\om,x) M_\om T_x\phi\,\hbox{d} x\,\hbox{d}\om.
\end{equation}
That is, we are using the time-frequency translates (i.e.,
translates and modulates) of $\phi$, $M_\om T_x\phi$, as a `basic
set' of functions by means of which one can express an arbitrary
$L^2$ function $f$.

However, the set $\{M_\om T_x\phi\hbox{:}\ \om,x\in\BR^d\}$ is
highly overcomplete, and in general one can expand an arbitrary
$f$ using only a suitable discrete subset of it. For example, if
$\phi$ is the characteristic function of the unit cube $[0,1]^d$,
the set $\{M_jT_k\phi\hbox{:}\ j,k\in\BZ^d\}$ is actually an
orthonormal basis for $L^2(\BR^d)$. (It is the basis one obtains
by tiling $\BR^d$ by cubes of unit side length with vertices in
the integer lattice and using the usual Fourier basis on each
cube.) Such discrete `basic sets' have an obvious advantage from a
computational point of view, as integrals must be approximated by
discrete sums for numerical work anyway. We are thus led to the
following definition.

Given $\phi\in L^2(\BR^d)$ and $\a,\b>0$, let
\begin{equation*}
\CG(\phi,\a,\b)=\blb M_{\a j}T_{\b k}\phi\hbox{:}\
j,k\in\BZ^d\brb.
\end{equation*}
$\CG(\phi,\a,\b)$ is called the \emph{Gabor system} determined by
$(\phi,\a,\b)$. (The name is in honor of the electrical engineer
D~Gabor, who suggested in his ground-breaking paper of 1946
\cite{gab} that $\CG(\hbox{e}^{-\pi t^2/\a^2},1/\a,\a)$ should be
a useful `basic set'. He was not entirely correct~--- see \PAR3.4
of \cite{fol1}~--- but the essential idea was still a good one.)
Evidently this concept can be~generalized: one can consider
$\{M_\om T_x\phi\hbox{:}\ (\om,x)\in \L\}$ where $\L$ is a more
general discrete subset of $\BR^{2d}$. We shall consider such
generalizations later, but we stick with the lattice
$\L=\a\BZ^d\times\b\BZ^d$ for now.

The first question that must be addressed is the following:
\emph{For which $\phi,$ $\a,$ and $\b$ does $\CG(\phi,\a,\b)$ span
$L^2(\BR^d)$}? (By `span' we mean that its finite linear span is
dense in $L^2(\BR^d)$.) This question, as it stands, is too broad
to admit a reasonable answer. In particular, if one is given
$\phi$, the set of $(\a,\b)$ for which $\CG(\phi,\a,\b)$ spans
$L^2$ depends strongly on $\phi$. For example, if $d=1$ and $\phi$
(resp.\ $\hat\phi$) is supported in an interval of length $l$, for
$\CG(\phi,\a,\b)$ to span $L^2$ it is obviously necessary that
$\b\le l$ (resp.\ $\a\le l$). On the other hand, no matter what
$\phi$ is, it is always necessary that the lattice
$\a\BZ^d\times\b\BZ^d$ should be sufficiently dense in
$\BR^d\times\BR^d$. Indeed, we have the following:

\begin{theor}[\!]
If $\a\b>1${\rm ,} there is no $\phi\in L^2$ such that
$\CG(\phi,\a,\b)$ spans $L^2(\BR^d)$.
\end{theor}

When $\a\b>1$ and $\a\b$ is rational, Daubechies \cite{daub1} has
shown how to produce, for any given $\phi$, an explicit $f\ne 0$
such that $f\perp\CG(\phi,\a,\b)$. (See also p.~107 of
\cite{daub2} for the very easy case where $\a=1$, $\b=2$.)
However, when $\a\b$ is irrational, a deeper argument is
necessary. The reason for this rather surprising distinction lies
in the structure of the group of operators $G_{\a,\b}$ generated
by the modulations and translations $M_{\a j}$ and $T_{\b k}$
($j,k\in\BZ^d$). In view of \eqref{3}, we have
\begin{equation}\label{5}
G_{\a,\b}=\blb \hbox{e}^{\tpi\a\b l}M_{\a j}T_{\b k}\hbox{:}\
j,k\in\BZ^d,\ l\in \BZ\brb.
\end{equation}
This group is a homomorphic image of the \emph{discrete Heisenberg
group} $\BH_d$ whose underlying set is $\BZ^d\times\BZ^d\times\BZ$
and whose group law is
\begin{equation}\label{6}
(j,k,l)\cdot(j',k',l')=(j+j',\,k+k',\,l+l'+k\cdot j').
\end{equation}
(This is often written with the roles of $j$ and $k$ switched.)
Indeed, the map
\begin{equation}\label{15}
\pi_{\a,\b}(j,k,l)=\hbox{e}^{-\tpi\a\b l}M_{\a j}T_{\b k}
\end{equation}
is easily seen to be a unitary representation of $\BH_d$ whose
image is $G_{\a,\b}$. When $\a\b$ is irrational, $\pi_{\a,\b}$ is
an isomorphism of groups. On the other hand, when $\a\b$ is
rational, say $\a\b=p/q$ in lowest terms, we have
\begin{equation}\label{25}
\text{ker}(\pi_{\a,\b})=qZ\equiv \blb(0,0,ql)\hbox{:}\
l\in\BZ\brb,
\end{equation}
and hence
\begin{equation*}
G_{\a,\b}\cong \BH_d/qZ\cong \blb(j,k,l)\hbox{:}\ j,k\in\BZ^d,\
l\in \BZ/q\BZ\brb,
\end{equation*}
the group law again being given by \eqref{6} with addition mod $q$
in the last coordinate. Now, observe that $Z=\{(0,0,l)\hbox{:}\
l\in\BZ\}$ is both the center and the commutator subgroup of
$\BH_d$. The group $\BH_d/qZ$ is `almost Abelian': its commutator
subgroup is finite, and it has normal Abelian subgroups of finite
index (for example, $\{(j,qk,l)\hbox{:}\ j,k\in\BZ^d,\
l\in\BZ/q\BZ\}$). As a result, the harmonic analysis of $\BH_d/qZ$
can be reduced to \emph{Abelian} harmonic analysis. But $\BH_d$
itself is a discrete group that is \emph{not} `almost Abelian,'
and consequently it is not even type~I. This means that its
harmonic analysis exhibits various pathologies and involves queer
beasts such as von Neumann factors of type~II (see Chapter~7 of
\cite{fol2} for a fuller explanation of these matters).

In \PAR2 we develop the basic properties of the von Neumann
algebra generated by $G_{\a,\b}$, and in \PAR3 we discuss several
proofs of Theorem~1.

We now turn to a different question. Suppose that
$\CG(\phi,\a,\b)$ does span $L^2$; can we use it in an efficient
way to expand an arbitrary $f\in L^2$? The condition that allows
everything to work smoothly is that $\CG(\phi,\a,\b)$ should be a
\emph{frame} for $L^2$. The notion of frame was introduced by
Duffin and Schaeffer in 1952, but it was not fully exploited until
a third of a century later, when Daubechies {\it et~al}
\cite{d-g-m} showed how handy frames could be in signal analysis.
We make a brief detour into abstract Hilbert space theory to
explain this idea; see \cite{daub2} or \cite{gro} for a fuller
discussion.

A bit of notation that will be employed throughout this paper: If
$I$ is a discrete set, we denote functions on $I$ by lower-case
boldface letters such as $\bc$, and the value of $\bc$ at $i\in I$
is denoted by $c_i$.

Let $\CH$ be a separable Hilbert space. A~countable set
$\{e_i\}_{i\in I}\sbs\CH$ is called a \emph{frame} for $\CH$ if
there exist $C_1,C_2>0$ such that for all $f\in\CH$ we have the
`frame inequalities'
\begin{equation}\label{22}
C_1\|f\|_2\le\sum_{i\in I}|\langle f,e_i\rangle|^2\le C_2\|f\|_2.
\end{equation}
The second inequality in \eqref{22} means that the linear map $A$
from $\CH$ to functions on $I$ defined by $(Af)_i=\langle
f,e_i\rangle$ is bounded from $\CH$ to $l^2(I)$. Its adjoint
$A^*\hbox{:}\ l^2(I)\to\CH$ is easily seen to be $A^*\bc=\sum
c_ie_i$, where the series converges unconditionally. The
composition $S=A^*A\hbox{:}\ \CH\to\CH$, given by
\begin{equation*}
Sf=\sum \langle f,e_i\rangle e_i,
\end{equation*}
is called the \emph{frame operator}. Since $\langle
Sf,f\rangle=\sum |\langle f,e_i\rangle |^2$, the frame
inequalities \eqref{22} are equivalent to the operator
inequalities $C_1I\le S\le C_2I$; in particular, $S$ is
invertible, and its inverse satisfies $C_2^\inv I\le S^\inv\le
C_1^\inv I$. Since
\begin{align*}
\sum |\langle f,S^\inv e_i\rangle|^2=\sum |\langle S^\inv
f,e_i\rangle|^2=\langle S(S^\inv f),S^\inv f\rangle = \langle
S^\inv f,f\rangle,
\end{align*}
these inequalities in turn imply that $\{S^\inv e_i\}_{i\in I}$ is
again a frame (with frame constants $C_2^\inv$ and $C_1^\inv$); it
is called the \emph{dual frame}. The two frames $\{e_i\}$ and
$\{S^\inv e_i\}$ can now be used together to produce expansion
formulas for a general $f\in\CH$ in terms of either frame:
\begin{align*}
f&=S(S^\inv f)=\sum \langle S^\inv f,e_i\rangle e_i=\sum \langle
f,S^\inv e_i\rangle e_i,\\[.5pc]
f&=S^\inv(Sf)=S^\inv\left(\sum \langle f,e_i\rangle e_i\right) =
\sum \langle f,e_i\rangle S^\inv e_i.
\end{align*}

We now return to the Gabor system $\CG(\phi,\a,\b)$. If this
system is a frame, an easy calculation shows that the frame
operator $S$ commutes with every $M_{\a j}$ and $T_{\b k}$, so the
dual frame is again a Gabor system, namely,
$\CG(S^\inv\phi,\a,\b)$; we call $S^\inv\phi$ the \emph{dual
window} to\break $\phi$.

Incidentally, Theorem~1 implies that \emph{if there exists
$\phi\in L^2(\BR^d)$ such that $\CG(\phi,\a,\b)$ is a frame{\rm ,}
then $\a\b\le 1$}, but this result is easier to prove than
Theorem~1 itself (see p.~108 of \cite{daub2} or Corollary~7.5.1 of
\cite{gro}).

In the applications of Gabor systems, one is generally interested
in using windows with good time-frequency localization, that is,
windows $\phi$ for which both $\phi$ and $\hat\phi$ have
reasonably rapid decay at infinity, or~--- what is more or less
the same thing~--- that the windowed Fourier transform
$V_\g\phi(\om,x)$, for some nice fixed window $\g$, has rapid
decay in both $\om$ and $x$. Experience has shown that a good
quantitative measure of this decay is given by the norms that
characterize the so-called \emph{modulation spaces} $M_v^1$, which
are defined as follows.

A \emph{subexponential weight} on $\BR^{2d}$ is a function
$v\hbox{:}\ \BR^{2d}\to[0,\inf)$ of the form
$v(\x)=\hbox{e}^{\s(\r(\x))}$ where $\r$ is a seminorm on
$\BR^{2d}$ and $\s$ is a nonnegative concave function on
$[0,\inf)$ such that $\s(0)=0$ and $\lim_{r\to\inf}\s(r)/r=0$.
(Examples: $v_1(\om,x)=(1+|\om|+|x|)^a$, $v_2(\om,x)=(1+|\om|)^a$,
and $v_3(\om,x)=\hbox{e}^{|x|^b}$, where $a>0$ and $0<b<1$.) Such
a weight is always submultiplicative: $v(\x_1+\x_2)\le
v(\x_1)v(\x_2)$. Given a subexponential weight $v$, the
\emph{modulation space} $M^1_v$ is defined as
\begin{equation}\label{13}
M_v^1=\left\{ f\in L^2(\BR^d)\hbox{:}\ \|f\|_{1,v}=\int_{\BR^{2d}}
|V_\g f(\x)|v(\x)\,\hbox{d}\x<\inf\right\},
\end{equation}
where $\g$ is some fixed Schwartz-class window; $M_v^1$ turns out
to be independent of the choice of $\g$, as is the equivalence
class of the norm $\|f\|_{1,v}$. Evidently the notion of
modulation space can be generalized~--- for example, by using an
$L^p$ norm rather than the $L^1$ norm~--- but $M^1_v$ will suffice
for our purposes. For a detailed treatment of modulation spaces we
refer to Chapter~11 of \cite{gro} and \cite{f-g1}. Incidentally,
$M^1_1$ (i.e., $M^1_v$ where $v\equiv 1$) is the Feichtinger
algebra, often denoted by $S_0(\BR^d)$ (see \cite{f-g1}).

Now, suppose $\phi$ is a window such that $\CG(\phi,\a,\b)$ is a
frame for $L^2(\BR^d)$. If $\phi$ has good time-frequency
localization in the sense that $\phi\in M_v^1$ for some suitable
weight $v$, it is obviously desirable that the dual window
$S^\inv\phi$ should also belong to $M_v^1$. That this is indeed
the case is the second major result we wish to discuss.

\begin{theor}[\!]
Suppose $v$ is a subexponential weight on $\BR^{2d}$ and
$\a,\b>0$. If $\phi\in M_v^1$ is a window such that
$\CG(\phi,\a,\b)$ is a frame{\rm ,} then the frame operator $S$
maps $M_v^1$ bijectively onto itself. In particular{\rm ,}
$S^\inv\phi\in M_v^1$.
\end{theor}

This theorem was first proved by Feichtinger and Gr\"ochenig
\cite{f-g} in the case where $\a\b$ is rational. The general case
is a more recent result of Gr\"ochenig and Leinert \cite{g-l}; its
proof involves some abstract machinery that was not needed for the
rational case. The underlying reason for this dichotomy is the
same as in Theorem~1: one needs a result about a noncommutative
convolution on $\BZ^{2d}$ that is closely related to the group
$G_{\a,b}$ defined in \eqref{5}. When $\a\b$ is rational, these
results can be obtained by Abelian harmonic analysis, but the
general case requires a different approach.

In more detail, the crucial ingredient for the proof of Theorem~2
is a noncommutative analogue of a classic result of Norbert
Wiener. Wiener's theorem is usually stated as follows: \emph{If
$f$ is a continuous{\rm ,} nonvanishing{\rm ,} periodic function
on $\BR$ whose Fourier series is absolutely convergent{\rm ,} then
the Fourier series of $1/f$ is also absolutely convergent.}
However, by passing from $f$ to its sequence $\bc$ of Fourier
coefficients, this result can also be stated as a theorem about
the convolution algebra $l^1(\BZ)$. Indeed, taking into account
the fact that the nonvanishing of $f$ is equivalent to the
invertibility of the operator $g\mapsto fg$ on $L^2(\BR/\BZ)$, and
hence of the operator $\ba\mapsto\bc*\ba$ on $l^2(\BZ)$, Wiener's
theorem can be restated as follows: \emph{Suppose $\bc\in
l^1(\BZ)$ and the map $\ba\mapsto\bc*\ba$ is invertible as an
operator on $l^2(\BZ)$. Then $\bc$ is invertible in the
convolution algebra $l^1(\BZ)$.}

The setting for the noncommutative analogue of Wiener's theorem is
as follows. Given a real number $\g$, we define the operation of
\emph{$\g$-twisted convolution}, denoted by $\nat_\g$, on
$l^1(\BZ^d\times\BZ^d)$ (or $l^1(\BZ^{2d})$ for short) by
\begin{equation}\label{7}
(\ba\nat_\g
\bb)_{jk}=\sum_{l,m}a_{lm}b_{(j-l)(k-m)}\hbox{e}^{-\tpi\g(j-l)\cdot
m}.
\end{equation}
We also define an involution $\ba\mapsto \ba^{*_\g}$ on
$l^1(\BZ^{2d})$ by
\begin{equation}\label{8}
(\ba^{*_\g})_{jk}=\bar {a_{(-j)(-k)}}\hbox{e}^{-\tpi\g j\cdot k}.
\end{equation}
(Note: these definitions differ from the ones in \cite{g-l} by the
minus signs in the exponents. This is to compensate for the fact
that in \cite{g-l}, time-frequency shifts are written as $T_{\a
j}M_{\b k}$ rather than $M_{\a j}T_{\b k}$.) $l^1(\BZ^{2d})$ is a
Banach $*$-algebra with product $\nat_\g$ and involution $*_\g$.
We shall denote this algebra by $\CA_\g$:
\begin{equation*}
\CA_\g=\big( l^1(\BZ^{2d}), \nat_\g, *_\g\big).
\end{equation*}
Moreover, the obvious analogue of Young's inequality holds: if
$\ba\in l^1(\BZ^{2d})$, the operator
\begin{equation}\label{9}
L_{\bf a}(\bb)=\ba\nat_\g\bb
\end{equation}
is bounded on every $l^p(\BZ^{2d})$ with norm at most $\|\ba\|_1$.
The analogue of Wiener's theorem is as follows:

\begin{theor}[\!]
If $\ba\in \CA_\g${\rm ,} then the spectrum of the operator
$L_{\bf a}$ as an operator on $l^2(\BZ^{2d})$ is equal to its
spectrum as an operator on $l^1(\BZ^{2d})$. In particular{\rm ,}
if $L_{\bf a}$ is invertible on $l^2(\BZ^{2d})${\rm ,} then $\ba$
is invertible in the algebra $\CA_\g$.
\end{theor}

Section~4 is devoted to a sketch of the proofs of Theorems~2 and
3.

\section{von Neumann algebras generated by translations and
modulations}

We review some notation and terminology. If $\CH$ is a Hilbert
space, $\CB(\CH)$ denotes the algebra of bounded linear operators
on $\CH$. A~\emph{von Neumann algebra} on $\CH$ is a weakly closed
$*$-subalgebra of $\CB(\CH)$. If $\CM\sbs\CB(\CH)$, its
\emph{commutant} $\CM'$ is the von Neumann algebra of all
$B\in\CB(\CH)$ that commute with all $A\in\CM$. A~fundamental
theorem of von Neumann states that if $\CM$ is itself a von
Neumann algebra, then $(\CM')'=\CM$. Takesaki \cite{tak2} is a
good general reference for the von Neumann algebra theory needed
here. Nelson \cite{nel} has given a particularly nice proof of the
theorem just quoted.

Given $\a,\b>0$, let $\CM_{\a,\b}$ be the von Neumann algebra on
$L^2(\BR^d)$ generated by the operators $M_{\b j}$ and $T_{\a k}$
($j,k,\in\BZ^d$), that is, the von Neumann algebra generated by
the group $G_{\a,\b}$. The first fundamental fact about
$\CM_{\a,\b}$ is the following.

\begin{propo}\label{p1}$\left.\right.$\vspace{.5pc}

\noindent $\CM'_{\a,\b}=\CM_{1/\b,1/\a}$.
\end{propo}

This is a special case of a theorem of Takesaki \cite{tak1}, of
which there are several proofs in the literature (see \cite{rie}).
It is obvious that $\CM_{1/\b,1/\a}\sbs\CM'_{\a,\b}$ because
$M_\om$ commutes with $T_x$ precisely when $\om\cdot x\in\BZ$.
A~simple proof of the reverse inclusion can be found in
Appendix~6.1 of \cite{d-l-l}. Given $S\in\CM'_{\a,\b}$ and
$T\in\CM'_{1/\b,1/\a}$, one shows by obtaining explicit
representations for $S$ and $T$ that $ST=TS$; hence
$\CM'_{\a,\b}\sbs(\CM'_{1/\b,1/\a})'=\CM_{1/\b,1/\a}$.

Next, we recall that if $\CA$ is a $*$-subalgebra of $\CB(\CH)$, a
\emph{faithful trace} on $\CA$ is a linear functional $\t$ on
$\CA$ such that $\t(AB)=\t(BA)$ for $A,B\in\CA$ and $\t(A^*A)>0$
for all nonzero $A\in\CA$. If $\t(I)=1$, $\t$ is said to be
\emph{normalized}.

Let $\CM^0_{\a,\b}$ be the finite linear span of the operators
$M_{\a j}T_{\b k}$ ($j,k\in\BZ^d$), and define $\t$ on
$\CM^0_{\a,\b}$ by
\begin{equation}\label{14}
\t\left(\sum c_{jk}M_{\a j}T_{\b k}\right)=c_{00}.
\end{equation}
It is easily verified that $\t$ is a normalized faithful trace on
$\CM^0_{\a,\b}$. Moreover, suppose $R_1,\ldots,R_N$ are
rectangular solids (products of intervals) in $\BR^d$ whose
interiors are disjoint, whose side lengths are all at most
$\min(1/\a,\b)$, and whose union is the cube $[0,1/\a]^d$. If
$\chi_n$ is the characteristic function of $R_n$, then
\begin{equation}\label{11}
\t(A)=\a^d\sum_1^N\langle A\chi_n,\chi_n\rangle,\quad
A\in\CM^0_{\a,\b}.
\end{equation}
Indeed, it suffices to verify \eqref{11} when $A=M_{\a j}T_{\b
k}$. Since the side lengths of the $R_n$ are at most $\b$, their
translates by amounts $\b k$ ($k\ne0$) are disjoint, so $\langle
M_{\a j}T_{\b k}\chi_n,\chi_n\rangle=0$ unless $k=0$, in which
case (since $\bigcup_1^N R_j=[0,1/\a]^d$)
\begin{equation*}
\a^d\sum\langle M_{\a
j}\chi_n,\chi_n\rangle=\a^d\int_{[0,1/\a]^d}\hbox{e}^{\tpi \a
j\cdot t}\,\hbox{d}t=\d_{j0}.
\end{equation*}

Equation \eqref{11} shows that $\t$ extends uniquely to a
normalized faithful trace on $\CM_{\a,\b}$ that is continuous in
the weak operator topology, so that $\CM_{\a,\b}$ is a
\emph{finite} von Neumann algebra.

We remark that if $\a\b$ is rational, say $\a\b=p/q$, then the
center of $\CM_{\a,\b}$ (that is,
$\CM_{\a,\b}\cap\CM_{1/\b,1/\a}$) is large: it contains all
operators $M_{\a qj}T_{\b qk}$ with $j,k\in\BZ^d$. However, if
$\a\b$ is irrational, then the center of $\CM_{\a,\b}$ is trivial;
that is, $\CM_{\a,\b}$ is a \emph{factor}. (The ideas used in
\cite{d-l-l} to prove Proposition~1, as sketched above, easily
yield a proof of this.) Since $\CM_{\a,\b}$ has a faithful trace,
it is actually a factor of type ${\rm II}_1$.

We need one further ingredient. Given $\phi\in L^2(\BR^d)$, let
$A_\phi$ be the map from $L^2(\BR^d)$ to the space of functions on
$\BZ^{2d}$ given by
\begin{equation}\label{16}
(A_\phi f)_{jk}=\langle f,M_{\a j}T_{\b k}\phi\rangle.
\end{equation}
(We encountered such maps earlier in the discussion of frames.)
If $A_\phi$ is bounded from $L^2(\BR^d)$ to $l^2(\BZ^{2d})$, its
adjoint is given by
\begin{equation*}
A_\phi^*\bc=\sum c_{jk}M_{\a j}T_{\b k}\phi,
\end{equation*}
where the series is unconditionally convergent.\pagebreak

Suppose $\phi,\psi\in L^2(\BR^d)$ and $A_\phi,A_\psi$ are bounded
from $L^2$ to $l^2$. We can then consider the `generalized frame
operator' $S_{\psi,\phi}=A_\psi^* A_\phi$:
\begin{align}\label{17}
S_{\psi,\phi} f=\sum \langle f,M_{\a j}T_{\b k}\phi\rangle M_{\a
j}T_{\b k}\psi.
\end{align}
An easy calculation that we leave to the reader shows that
$S_{\psi,\phi}$ commutes with every $M_{\a j}$ and $T_{\b k}$, so
by Proposition~1, $S_{\psi,\phi}\in \CM_{1/\b,1/\a}$. The explicit
expansion of $S_{\psi,\phi}$ in terms of the operators
$M_{j/\b}T_{k/\a}$ is quite pretty; it is known as the
\emph{Janssen representation}:
\begin{align}\label{12}
S_{\psi,\phi}=(\a\b)^{-d}\sum_{j,k} \langle
\psi,M_{j/\b}T_{k/\a}\phi\rangle M_{j/\b}T_{k/\a}.
\end{align}
Actually, without additional hypotheses on $\psi$ and $\phi$ the
convergence of the series on the right is questionable, but it is
sufficient for $\psi$ and $\phi$ to belong to the Feichtinger
algebra $M^1_1$ defined by \eqref{13} with $v\equiv 1$.

\begin{propo}\label{p2}$\left.\right.$\vspace{.5pc}

\noindent If $\phi$ and $\psi$ are in $M^1_1${\rm ,} then $A_\phi$
and $A_\psi$ are bounded from $L^2(\BR^d)$ to $l^2(\BZ^{2d})$.
Moreover{\rm ,} the Janssen representation \eqref{12} for
$S_{\psi,\phi}$ is valid{\rm ,} and the series on the right
converges absolutely in the operator norm.
\end{propo}

For the proof, see Theorem~7.2.1 and Proposition~12.1.11 of
\cite{gro}.

The algebra $\CM_{1/\b,1/\a}$ has a normalized faithful trace just
like $\CM_{\a,\b}$, which we denote by $\t'$. The formula for the
trace of $S_{\psi,\phi}$ is very simple.

\begin{propo}\label{p3}$\left.\right.$\vspace{.5pc}

\noindent If $A_\phi$ and $A_\psi$ are bounded from $L^2(\BR^d)$
to $l^2(\BZ^{2d})${\rm ,} then
\begin{equation*}
\t'(S_{\psi,\phi})=(\a\b)^{-d}\langle \psi,\phi\rangle.
\end{equation*}
\end{propo}

This is an immediate corollary of \eqref{14} (with $1/\b$, $1/\a$
in place of $\a$, $\b$) and the Janssen representation if
$\phi,\psi\in M_1^1$. The general case is not hard to prove from
formula \eqref{11} (see \cite{d-l-l}).

\section{Time-frequency density of complete Gabor systems}

We recall that if $\CM$ is a von Neumann algebra on $\CH$, a
vector $v\in\CH$ is called \emph{cyclic} for $\CM$ if
$\{Av\hbox{:}\ A\in\CM\}$ is dense in $\CH$, and $\CM$ is called
\emph{cyclic} if it has a cyclic vector. The Gabor system
$\CG(\phi,\a,\b)$ spans $L^2(\BR^d)$ precisely when $\phi$ is a
cyclic vector for $\CM_{\a,\b}$, so Theorem~1 can be restated as
follows:
\begin{equation}\label{10}
\text{\emph{If $\CM_{\a,\b}$ is cyclic, then
$\a\b\le 1$.}}
\end{equation}

The first proof of Theorem~1, in a sense, appeared before the
question was even posed. That is, \eqref{10} is a consequence of a
theorem of Rieffel \cite{rie} concerning the `coupling function'
for (a generalization of) $\CM_{\a,\b}$ and its commutant. The
connection with Rieffel's theorem was first pointed out in
Daubechies \cite{daub1}. However, Rieffel's paper \cite{rie} is
quite technical, and for those who are not specialists in operator
algebras (including the present writer) the arguments in it, and
the notion of `coupling function' itself, are hard to grasp. We
shall not attempt to describe them further.

At about the same time as \cite{daub1}, another proof of
Theorem~1, for the case $d=1$, was given by Baggett \cite{bag}.
Baggett's argument is of interest because it explicitly develops
the connection of the problem with the representation theory of
the discrete Heisenberg group $\BH_1$. (It can be generalized to
the $d$-dimensional case.)

We recall that the group $G_{\a,\b}$ generated by the operators
$M_{\a j}$ and $T_{\b k}$ is the image of $\BH_1$ under the
unitary representation $\pi_{\a,\b}$ defined by \eqref{15}.
Baggett relates $\pi_{\a,\b}$ to some other representations of
$\BH_1$. First, let $H$ be the Abelian subgroup of $\BH_1$
consisting of elements of the form $(j,0,l)$ ($j,l\in\BZ$), and
for $\g,\d\in\BR$ let $\s_{\g,\d}$ be the representation of
$\BH_1$ induced from the character
$\chi_{\g,\d}(j,0,l)=\hbox{e}^{\tpi(\d j-\g l)}$ of $H$.
$\s_{\g,\d}$ acts on $l^2(\BZ)$ by
\begin{equation*}
[\s_{\g,\d}(j,k,l)\bc]_n=\hbox{e}^{\tpi(\d j-\g l+n\g j)}c_{n-k}.
\end{equation*}
Next, for $\y>0$ let $\S_{\g,\y}$ be the representation of $\BH_1$
defined as the direct integral
\begin{equation*}
\S_{\g,\y}=\int_{[0,\y)}^\oplus \s_{\g,\d}\,\hbox{d}\d,
\end{equation*}
which acts on $L^2([0,\y)\times\BZ)$ in the obvious way, and let
$\CN_{\g,\y}$ be the von Neumann algebra on $L^2([0,\y)\times\BZ)$
generated by the operators $\S_{\g,\y}(h)$, $h\in\BH_1$. Baggett
establishes the following facts:

\begin{enumerate}
\renewcommand{\labelenumi}{(\roman{enumi})}
\leftskip .4pc
\item $\pi_{\a,\b}$ is unitarily equivalent to $\S_{\g,\g}$ where
$\g=\a\b$.
\item $\S_{\g,\y}$ is unitarily equivalent to $\S_{\g',\y'}$ if and only
if $\g=\g'$ and $\y=\y'$.
\item $\CN_{\g,\y}$ is isomorphic to $\CN_{\g,\g}$ whenever $\y\ge\g$.
\item $\S_{\g,1}$ is unitarily equivalent to the representation of
$\BH_1$ induced from the central character
$\chi_\g(0,0,l)=\hbox{e}^{-\tpi\g l}$.
\end{enumerate}

The proof of Theorem~1, in the form \eqref{10}, now proceeds by
contradiction. Indeed, suppose that $\g=\a\b>1$ but $\CM_{\a,\b}$
is cyclic. Then $\CM'_{\a,\b}=\CM_{1/\b,1/\a}$ is also cyclic
since $1/\a\b<1$. (If $1/\a\b<1$, it is easy to construct $\phi\in
L^2(\BR^d)$ such that $\CG(\phi,1/\b,1/\a)$ is a frame for
$L^2(\BR^d)$; see \cite{d-g-m}.) Hence, by (i), $\CN_{\g,\g}$ and
its commutant are both cyclic. On the other hand, from (iv) and
the theory of induced representations, $\CN_{\g,1}$ and its
commutant are also both cyclic. But then, by well-known facts
about von Neumann algebras that can be found in \cite{tak2}, the
fact (iii) that $\CN_{\g,\g}$ and $\CN_{\g,1}$ are isomorphic
implies that they are actually unitarily equivalent, which
contradicts (ii).

A few years after Daubechies \cite{daub1} and Baggett \cite{bag},
Daubechies {\it et~al} \cite{d-l-l} found another proof of
Theorem~1 that uses the trace $\t'$ on $\CM_{1/\b,1/\a}$ and the
operators $A_\phi$ and $S_{\psi,\phi}$ defined by \eqref{16} and
\eqref{17} in a very efficient way. (Like \cite{bag}, this paper
deals explicitly only with the case $d=1$, but the generalization
to arbitrary $d$ is entirely straightforward.)

Suppose that $\psi\in L^2(\BR^d)$ is a cyclic vector for
$\CM_{\a,\b}$. If $A_\psi$ is not bounded from $L^2(\BR^d)$ to
$l^2(\BZ^{2d})$, we regard it as an unbounded linear map with
domain $D(A_\psi)=\{f\in L^2\hbox{:}\ \sum|\langle f,M_{\a j}T_{\b
k}\psi\rangle|^2<\inf\}$; as such, it is always closed and densely
defined. It follows, by a theorem of von Neumann, that
$S_{\psi,\psi}=A_\psi^*A_\psi$ is a positive self-adjoint operator
on $L^2(\BR^d)$. As in the bounded case, it commutes with the
$M_{\a j}$ and $T_{\b k}$, and so does $(\e I+S_{\psi,\psi})^\inv$
for any $\e>0$. It follows that if we set
\begin{equation*}
\phi=(\e I+S_{\psi,\psi})^\inv\psi,
\end{equation*}
then for any $f\in L^2(\BR^d)$,
\begin{align*}
(A_\phi f)_{jk}&=\langle f,M_{\a j}T_{\b k}(\e
I+S_{\psi,\psi})^\inv\psi\rangle\\[.5pc]
& =\langle (\e I+S_{\psi,\psi})^\inv f,M_{\a j}
T_{\b k}\psi\rangle=(A_\psi(\e I+S_{\psi,\psi})^\inv f)_{jk},
\end{align*}
that is,
\begin{equation*}
A_\phi=A_\psi(\e I+S_{\psi,\psi})^\inv,
\end{equation*}
and hence
\begin{equation}\label{19}
S_{\psi,\psi}(\e I+S_{\psi,\psi})^\inv=A_\psi^*A_\phi.
\end{equation}
Next, Daubechies {\it et~al} \cite{d-l-l} showed that the formula
in Proposition~3 for the trace of $S_{\psi,\phi}=A_\psi^*A_\phi$
remains valid provided only $S_{\psi,\phi}$ is a bounded operator,
even if $A_\psi$ or $A_\phi$ is not. By \eqref{19}, that is the
case for the $\psi$ and $\phi$ under consideration here, so
\begin{equation*}
(\a\b)^d\t'(S_{\psi,\psi}(\e
I+S_{\psi,\psi})^\inv)=(\a\b)^d\t'(A_\psi^* A_\phi)=\langle
\psi,\phi\rangle.
\end{equation*}
Now, by the spectral functional calculus, $S_{\psi,\psi}(\e
I+S_{\psi,\psi})^\inv$ converges strongly as $\e\to0$ to the
orthogonal projection onto the closure of the range of
$S_{\psi,\psi}$. But since $\psi$ is cyclic for $\CM_{\a,\b}$,
this range is dense in $L^2$, so $S_{\psi,\psi}(\e
I+S_{\psi,\psi})^\inv$ converges strongly to $I$. The formula
analogous to \eqref{11} for $\t'$ shows that $\t'$ is continuous
with respect to strong convergence, so
\begin{equation*}
(\a\b)^d=(\a\b)^d\t'(I)=\langle \psi,\phi\rangle.
\end{equation*}
On the other hand, we have
\begin{equation*}
\langle \psi,\phi\rangle=\langle (\e I+A_\psi^*
A_\psi)\phi,\phi\rangle=\e\|\phi\|^2+\|A_\psi\phi\|^2\ge\|A_\psi\phi\|^2,
\end{equation*}
and if $\be\in l^2(\BZ^{2d})$ is defined by $e_{jk}=\d_{j0}\d_{k0}$,
\begin{equation*}
\langle \psi,\phi\rangle=\langle A_\psi^*\be,\phi\rangle=\langle
\be,A_\psi\phi\rangle\le\|A_\psi\phi\|.
\end{equation*}
(No absolute values are needed since $\langle
\psi,\phi\rangle=(\a\b)^d>0$.) These two inequalities imply that
$\langle \psi,\phi\rangle^2\le\langle \psi,\phi\rangle$ and hence
$(\a\b)^d=\langle \psi,\phi\rangle\le 1$. Thus the proof is
complete.

At about the same time as Daubechies {\it et~al} \cite{d-l-l},
Ramanathan and Steger \cite{r-s} found a more elementary argument
to prove, and indeed generalize, Theorem~1. They avoid von Neumann
algebras but use an idea from a different branch of harmonic
analysis: the notion of asymptotic density of a discrete set first
exploited by Beurling in his work on balayage for the Fourier
transform \cite{beur}.

To wit, let $\L$ be a discrete subset of
$\BR^{2d}=\BR^d\times\BR^d$, and for $\phi\in L^2(\BR^d)$, let
$\CG(\phi,\L)$ be the corresponding Gabor system:
\begin{equation*}
\CG(\phi,\L)=\blb M_\om T_x \phi\hbox{:}\ (\om,x)\in\L\brb.
\end{equation*}
We are interested in the question of whether $\CG(\phi,\L)$ spans
$L^2$. For $r>0$ and $(\y,a)\in\BR^{2d}$, let $B_r(\y,a)$ be the
ball of radius $r$ about $(\y,a)$, and let
$v(r)=(2\pi^d/d!)r^{2d}$ be its volume. Define
\begin{equation*}
\n^-(r)=\min_{(\y,a)\in\BR^{2d}}\mbox{card}(\L\cap B_r(\y,a)),
\end{equation*}
and define the \emph{lower density} of $\L$ to be
\begin{equation*}
D^-(\L)=\liminf_{r\to\inf} \frac{\n^-(r)}{v(r)}.
\end{equation*}
(There is a corresponding notion of upper density, but we shall
not need it.)

We shall say that $\CG(\phi,\L)$ has the \emph{homogeneous
approximation property} if for every $f\in L^2(\BR^d)$ and $\e>0$
there is an $R>0$ such that for every $(\y,a)\in \BR^{2d}$ there
is a finite linear combination $h$ of the $M_\om T_x \phi$'s with
$(\om,x)\in\L\cap B_R(\y,a)$ such that $\|h-M_\y T_a f\|_2<\e$.
(This definition is quite a mouthful. The homogeneous
approximation property implies that $\CG(\phi,\L)$ spans
$L^2(\BR^d)$, as one sees simply by taking $(\y,a)=(0,0)$. But it
is stronger: it means not only that each $f\in L^2$ can be
approximated by finite linear combinations of $M_\om T_x\phi$'s
with $(\om,x)\in\L$, but that translates and modulates of $f$ by
arbitrary amounts $(\y,a)$ can be uniformly approximated by linear
combinations of $M_\om T_x\phi$'s with $(\om,x)\in\L$ not too far
from $(\y,a)$.)

The main result of Ramanathan and Steger \cite{r-s} is the
following:

\begin{theor}[\!]
If $\CG(\phi,\L)$ has the homogeneous approximation property{\rm
,} then $D^-(\L)\ge1$.
\end{theor}

\noindent Let $\chi$ be the characteristic function of the unit
cube. The idea of the proof is to compare $\CG(\phi,\L)$ with
$\CG(\chi,\BZ^{2d})$, which is an orthonormal basis for
$L^2(\BR^d)$. For $(\y,a)\in\BR^{2d}$ and $r>0$, let
\begin{align*}
V_r(\y,a)&=\mbox{linear span of }\blb M_jT_k\chi\hbox{:}\
(j,k)\in\BZ^{2d}\cap B_r(\y,a)\brb,\\[.4pc]
W_r(\y,a)&=\mbox{linear span of }\blb M_\om T_x \phi\hbox{:}\
(\om,x)\in\L\cap B_r(\y,a)\brb.
\end{align*}
Next, given $R>0$, let $T$ be the restriction to the
finite-dimensional space $V_r(\y,a)$ of the operator
$P_{V_r(\y,a)}\circ P_{W_{r+R}(\y,a)}$ ($P_X$ = orthogonal
projection onto $X$). Since $T$ is a composition of projections,
its eigenvalues lie in $[0,1]$, so its trace is dominated by its
rank. Thus,
\begin{equation*}
\mbox{tr}(T)\le \dim W_{r+R}(\y,a) = \mbox{card}(\L\cap B_{r+R}(\y,a)).
\end{equation*}
On the other hand, the homogeneous approximation property implies
that for any $\e>0$ we can find $R>0$ such that
$\|T(M_jT_k\chi)-M_jT_k\chi\|_2<\e$ for all $(j,k)\in\BZ^{2d}\cap
B_r(\y,a)$, and hence
\begin{align*}
\mbox{tr}(T)=\sum_{\BZ^{2d}\cap B_r(\y,a)}\langle
TM_jT_k\chi,M_jT_k\chi\rangle\ge (1-\e)\mbox{card}(\BZ^{2d}\cap
B_r(\y,a)).
\end{align*}
Therefore,
\begin{equation*}
(1-\e)\frac{\mbox{card}(\BZ^{2d}\cap B_r(\y,a))}{v(r)} \le
\frac{\mbox{card}(\L\cap B_{r+R}(\y,a))}{v(r+R)}
\frac{v(r+R)}{v(r)}.
\end{equation*}
Pick sequences $r_n\to\inf$ and $(\y_n,a_n)\in\BR^{2d}$ such that
$\mbox{card}(\L\cap B_{r_n+R}(\y_n,a_n))/v(r_n+R)\to D^-(\L)$. We
have $\mbox{card}(\BZ^{2d}\cap B_{r_n}(\y_n,a_n))/v(r_n)\to 1$ and
$v(r_n+R)/v(r_n)=(r_n+R)^d/r_n^d\to 1$, and hence
\begin{equation*}
1-\e\le D^-(\L).
\end{equation*}
Since $\e$ is arbitrary, the proof is complete.\pagebreak

Ramanathan and Steger further show that \emph{if $\L$ is a lattice
and $\CG(\phi,\L)$ spans $L^2(\BR^d)$, then $\CG(\phi,\L)$ has the
homogeneous approximation property.} (The proof is
straightforward.) Since $D^-(\L)=1/{\rm vol}(\BR^{2d}/\L)$ in this
case, we have the following generalization of Theorem~1.

\begin{theor}[\!]
If $\L$ is a lattice in $\BR^{2d}$ such that ${\rm
vol}(\BR^{2d}/\L)>1${\rm ,} there is no $\phi$ in $L^2(\BR^d)$
such that $\CG(\phi,\L)$ spans $L^2(\BR^d)$.
\end{theor}

More recently, Gabardo and Han \cite{g-h} have extended Theorem~1
to a very general setting using the same circle of ideas as
Daubechies {\it et~al} \cite{d-l-l}, but arranged rather
differently. They define a \emph{group-like unitary system} to be
a countable collection $\CU$ of unitary operators on a separable
Hilbert space $\CH$ with the property that the group $G(\CU)$
generated by $\CU$ lies in $\{\l U\hbox{:}\ U\in\CU, \l\in\BC\}$.
Thus, $\{M_{\a j}T_{\b k}\hbox{:}\ j,k\in\BZ^d\}$ is a group-like
unitary system on $L^2(\BR^d)$, and any projective representation
of a countable group gives rise to a group-like unitary system.

Given a group-like unitary system $\CU$ on $\CH$ and $v\in\CH$,
one has the linear map $A_v$ from $\CH$ to functions on $\CU$
given by $(A_vw)_U=\langle w,Uv\rangle$; Gabardo and Han assume
that the set of $v\in\CH$ such that $A_v$ is bounded from $\CH$ to
$l^2(\CU)$ is dense in $\CH$. Under this condition, they show that
$\CH$ is the orthogonal direct sum of subspaces $\CH_i$ ($i\in
I$), each of which is $\CU$-invariant, and each of which possesses
a vector $v_i$ such that $\|w\|^2=\sum_{U\in {\CCU}}|\langle
w,Uv_i\rangle|^2$ for all $w\in\CH_i$ (that is, $\{Uv_i\hbox{:}\
U\in\CU\}$ is a `normalized tight frame' for $\CH_i$). They show
that $\sum_{i\in I} \|v_i\|^2$ depends only on $\CU$, and they
define the \emph{redundancy} of $\CU$ to be
\begin{equation*}
r(\CU)=\left(\sum \|v_i\|^2\right)^\inv.
\end{equation*}
They then show that:

\begin{enumerate}
\leftskip .2pc
\renewcommand{\labelenumi}{(\roman{enumi})}
\item $\CU$ has a cyclic vector (a vector $v\in\CH$ such that the
linear span of $\{Uv\hbox{:}\ U\in\CU\}$ is dense in $\CH$) if and
only if $r(\CU)\ge 1$.
\item If the index set $I$ is finite, then the commutant $\CU'$ is a
finite von Neumann algebra. The formula
$\tilde\t(A_w^*A_v)=\langle w,v\rangle$ (for any $v$ and $w$ such
that $A_v$ and $A_w$ are bounded from $\CH$ to $l^2(\CU)$)
determines a faithful trace on $\CU'$, and $r(\CU)=\tilde\t(I)$.
\end{enumerate}

In the case $\CU=\{M_{\a j}T_{\b k}\hbox{:}\ j,k\in\BZ^d\}$, one
can take the vectors $v_i$ to be the characteristic functions
$\chi_n$ in \eqref{11} and the subspaces $\CH_i$ to be the
$\CU$-invariant subspaces they generate. Then
$\CU'=\CM_{1/\b,1/\a}$, and by Proposition~3, Gabardo and Han's
trace $\tilde\t$ is $(\a\b)^d\t'$, so their results (i) and (ii)
imply Theorem~1. In fact, without much additional effort, they
yield Theorem~5 for the case of lattices of the form $A\BZ^d\times
B\BZ^d$ where $A,B\in GL(n,\BR)$.

Another generalization of Theorem~1 has been obtained by Bekka
\cite{bek}. Bekka considers an irreducible square-integrable
representation $\pi$ of a unimodular locally compact group $G$ on
a separable Hilbert space $\CH$, and a discrete subgroup $\G$ of
$G$ such that the volume of $G/\G$ is finite. He proves results
relating (i) the formal dimension of $\pi$ (i.e., the constant
$d_\pi$ such that $\|\xi\|^2\|\y\|^2=d_\pi\int_G|\langle
\x,\pi(g)\y\rangle|^2\,dg$), (ii) the volume of $G/\G$, and (iii)
the `center-valued von Neumann dimension' of $\CH$,
$\hbox{cdim}(\CH)$, as a VN($\G$)-module, where VN($\G$) is the
von Neumann algebra on $l^2(\G)$ generated by the left regular
representation. (We shall not attempt to describe $\rm{cdim}(\CH)$
other than to say that it is an element of the center of
VN($\G$).) Bekka's abstract version of Theorem~1 in this setting
is as follows:

\begin{theor}[\!]
If there is a vector $v\in\CH$ such that \hbox{$\{\pi(g)v{\rm :}\
g\in\G\}$} spans $\CH${\rm ,} then\break ${\rm cdim}(\CH)\le I$.
\end{theor}

When specialized to the case where $G$ is the (reduced) real
Heisenberg group and $\G$ is the discrete Heisenberg group, this
yields Theorem~5. (Note that Bekka normalizes his frequency
variables differently, with the result that some $2\pi$'s appear
in his formulas that are not in ours. But in fact each of his
$2\pi$'s should be $(2\pi)^d$.)

\section{Wiener's theorem and Gabor frames}

We now turn to the proofs of Theorems~2 and 3, following the
arguments of Gr\"ochenig and Leinert \cite{g-l}. We begin with
Theorem~3 and then obtain Theorem~2 from it.

Recall that we are concerned with the Banach $*$-algebra $\CA_\g$,
which is $l^1(\BZ^{2d})$ equipped with the product and involution
defined by \eqref{7} and \eqref{8}, and that for $\ba\in\CA_\g$,
$L_{\bf a}$ denotes the operator $L_{\bf a}(\bb)=\ba\nat_\g \bb$.
If $\CX$ is a space of functions on $\BZ^{2d}$ on which $L_{\bf
a}$ is a bounded operator, we denote the spectrum and spectral
radius of $L_{\bf a}$ on $\CX$ by $\s_{\CCX}(\ba)$ and
$\r_{\CCX}(\ba)$, respectively. That is,
\begin{align*}
\s_{\CCX}(\ba) &= \blb\l\in\BC\hbox{:}\ \l I-L_\ba\mbox{ is not
invertible on}\CX\brb,\\[.5pc]
\r_{\CCX}(\ba) &= \sup\blb |\l|\hbox{:}\ \l\in\s_{\CCX}(\ba)\brb =
\lim_{n\to\inf}\|L_\ba^n\|^{1/n}_{\CCB(\CCX)}.
\end{align*}
Here and in the sequel, $\CB(\CX)$ denotes the space of bounded
linear operators on $\CX$. Moreover, we abbreviate $l^p(\BZ^{2d})$
as $l^p$. Our goal is to prove that $\s_{l_1}(\ba)=\s_{l^2}(\ba)$
for all $\ba\in l^1=\CA_\g$.

Just as twisted convolution on $\BR^{2d}$ is closely related to
ordinary convolution on the~reduced Heisenberg group (the real
Heisenberg group modified so that its center is a circle rather
than a line; see pp.~25--26 of \cite{fol1}), so the twisted
convolution $\nat_\g$ on $\BZ^{2d}$ is closely related to ordinary
convolution on a modification of the discrete Heisenberg group
$\BH_d$, namely, the group $H_\g$ whose underlying set is
$\BZ^d\times\BZ^d\times\BT$ ($\BT$ = the group of complex numbers
of modulus one) and whose group law is
\begin{equation*}
(j,k,\z)(j',k',\z')=(j+j',\,k+k',\,\z\z'\hbox{e}^{-\tpi \g k\cdot
j'}).
\end{equation*}
(Note that $(j,k,l)\mapsto(j,k,\hbox{e}^{-\tpi\g l})$ is a
homomorphism from $\BH_d$ to $H_\g$. It is injective if $\g$ is
irrational, and its kernel is $qZ$ as in \eqref{25} if $\g=p/q$ in
lowest terms.) Indeed, $L^1(H_\g)$ is a Banach $*$-algebra under
convolution and the usual involution $f^*(\xi)=\bar{f(\xi^\inv)}$,
and the map $J\hbox{:}\ \CA_\g\to L^1(H_\g)$ defined by
\begin{equation}\label{20}
J(\ba)(j,k,\z)=\z^\inv a_{jk}
\end{equation}
is easily seen to be a $*$-isomorphism of $\CA_\g$ onto the
subalgebra of $L^1(H_\g)$ consisting of those functions $f$ that
satisfy $f(j,k,\z)=\z^\inv f(j,k,1)$.

With this in mind, the ingredients from abstract harmonic analysis
that are needed to prove Theorem~2 are as follows. In them we
employ an obvious modification of the notation introduced earlier:
if $f$ is an element of the Banach algebra $\CA$, $\s_{\CCA}(f)$
and $\r_{\CCA}(f)$ denote its spectrum and spectral radius. The
first two lemmas are theorems of Hulanicki \cite{hul} and Ludwig
\cite{lud}, respectively.

\begin{lem}\hskip -.4pc{\rm \cite{hul}.}\ \
Suppose $\CS$ is a $*$-subalgebra of a Banach $*$-algebra
$\CA${\rm ,} and there exists a faithful $*$-representation of
$\CA$ on a Hilbert space $\CH$ such that if $f\in\CS$ and $f=f^*$
then $\|\pi(f)\|_{\CCB(\CCH)}=\r_{\CCA}(f)$. Then for any
$f\in\CS$ with $f=f^*$ we have
$\s_{\CCA}(f)=\s_{\CCB(\CCH)}(\pi(f))$.
\end{lem}

Recall that a Banach $*$-algebra $\CA$ is called \emph{symmetric}
if $\s_{\CCA}(f)\sbs[0,\inf)$ for all $f\in \CA$ such that
$f=f^*$.

\begin{lem}\hskip -.4pc{\rm \cite{lud}.}\ \
If $G$ is a locally compact nilpotent group{\rm ,} then $L^1(G)$
is symmetric.
\end{lem}

Using the fact that every locally compact nilpotent group is
amenable, some standard facts about amenable groups, and Lemmas~1
and 2, it is not hard to deduce the next result (Theorem~2.8 of
\cite{g-l}).

\begin{lem}
Suppose $G$ is a locally compact nilpotent group and $f\in L^1(G)$
satisfies $f=f^*$. Then the convolution operator $C_f(g)=f*g$
satisfies $\s_{\CCB(L^1(G))}(C_f)=\s_{\CCB(L^2(G))}(C_f)$.
\end{lem}

We can now sketch the proof of Theorem~3. Suppose $\ba\in \CA_\g$.
If $\ba=\ba^{*_\g}$, we employ the map $J$ defined by \eqref{20}
to transfer the problem to $L^1(H_\g)$. Since $H_\g$ is nilpotent,
the conclusion $\s_{l^1}(\ba)=\s_{l^2}(\ba)$ follows easily from
Lemma~3. The result for general $\ba$ is now obtained by the
following simple device. If $L_\ba$ is invertible on $l^2$, then
so is $L_\bb$ where $\bb=a^{*_\g}\nat_\g\ba$ or $\bb=\ba\nat_\g
\ba^{*_\g}$, and these $\bb$'s satisfy $\bb=\bb^{*_\g}$, so
$\ba^{*_\g}\nat_\g \ba$ and $\ba\nat_\g\ba^{*_\g}$ are invertible
in $l^1$ by the result just proved. But then
$(\ba^{*_\g}\nat_\g\ba)^\inv\nat_\g\ba^{*_\g}$ and
$\ba^{*_\g}\nat_\g(\ba\nat_\g\ba^{*_\g})^\inv$ are left and right
inverses for $\ba$, respectively, so $\ba$ is invertible in $l^1$.
The desired conclusion $\s_{l^1}(\ba)=\s_{l^2}(\ba)$ follows by
applying this result to $\l\be-\ba$ where $\be$ is the identity in
the algebra $\CA_\g=l^1$ ($e_{jk}=\d_{j0}\d_{k0}$).

Finally, we show how Gr\"ochenig and Leinert \cite{g-l} deduce
Theorem~2 from Theorem~3. Let us suppose that $\CG(\phi,\a,\b)$ is
a frame for $L^2(\BR^d)$ and $\phi\in M^1_v$ where $v$ is a
subexponential weight function; we wish to show that the frame
operator $S=S_{\phi,\phi}$ (notation as in \eqref{17}) is
invertible on $M^1_v$. The link with Theorem~3 comes through the
following considerations.

\begin{lem}\label{l4}
If $\ba\in l^1(\BZ^{2d})${\rm ,} let
\begin{equation*}
\pi(\ba)=\sum a_{jk}M_{j/\b}T_{k/\a}.
\end{equation*}
Then $\pi$ is a $*$-representation of the algebra $\CA_\g$ where
$\g=1/\a\b$.
\end{lem}

Let
\begin{equation*}
\tilde v(j,k)=v(j/\b, k/\a),\quad j,k\in \BZ^d,
\end{equation*}
and let
\begin{equation*}
l^1_{\tilde v}=\blb \ba\in l^1(\BZ^{2d})\hbox{:}\ \sum
|a_{jk}|\tilde v(j,k)<\inf\brb.
\end{equation*}
Since $v$ is submultiplicative, it is easily verified that
$l^1_{\tilde v}$ is a $*$-subalgebra of the twisted convolution
algebra $\CA_\g$, for any $\g>0$.

\begin{lem}\label{l5}
If $\ba\in l^1_{\tilde v}${\rm ,} then $\pi(\ba)$ is bounded on
$M^1_v$.
\end{lem}

\begin{lem}\label{l6}
If $\ba\in l^1_{\tilde v}${\rm ,} then $\r_{l^1_{\tilde
v}}(\ba)=\r_{l^1}(\ba)$.
\end{lem}

The proofs of Lemmas~4, 5, and 6 are all quite easy. The next one
is a little deeper.

\begin{lem}\label{l7}
If $\ba\in l^1(\BZ^{2d})${\rm ,} then
$\|\pi(\ba)\|_{\CCB(L^2(\BR^d))}=\|L_\ba\|_{\CCB(l^2(\BZ^{2d}))}$.
\end{lem}

The idea of the proof is as follows. Let $C^*(l^1)$ be the $C^*$
subalgebra of $\CB(l^2(\BZ^{2d}))$ generated by the operators
$L_\ba$, $\ba\in l^1$, and let $C^*(\a,\b)$ be the $C^*$
subalgebra of $\CB(L^2(\BR^d))$ generated by the
$M_{j/\b}T_{k/\a}$, $j,k\in\BZ^d$. One shows that the
correspondence $L_\ba\mapsto \pi(\ba)$ extends to an injective
$*$-homomorphism from $C^*(l^1)$ to $C^*(\a,\b)$. But injective
$*$-homomorphisms of $C^*$ algebras are always isometries.

\begin{lem}\label{l8}
Suppose $\a,\b >0$. If $\ba\in L^1_{\tilde v}$ and $\pi(\ba)$ is
invertible on $L^2(\BR^d)${\rm ,} then $\ba$ is invertible in the
algebra $l^1_{\tilde v}\sbs\CA_\g${\rm ,} where $\g=1/\a\b$.
\end{lem}

To prove this, first suppose that $\ba=\ba^{*_\g}$. Then
\begin{equation*}
\|\pi(\ba)\|_{\CCB(L^2)}=\|L_\ba\|_{\CCB(l^2)}=\r_{l^2}(\ba)=\r_{l^1}(\ba)=\r_{l^1_{\tilde
v}}(\ba),
\end{equation*}
where the equalities are justified by Lemma~7, the fact that
$\ba=\ba^{*_\g}$, Theorem~3, and Lemma~6, respectively. Thus, we
can apply Hulanicki's theorem (Lemma~1) with $\CS=l^1_{\tilde v}$,
$\CA=\CA_\g$, $\pi$ as in Lemma~4, and $\CH=L^2(\BR^d)$ to
conclude that $\s_{l^1_{\tilde v}}(\ba)=\s_{\CCB(L^2)}(\pi(\ba))$.
Since $0\notin\s_{\CCB(L^2)}(\pi(\ba))$ by assumption, $\ba$ is
invertible in $l^1_{\tilde v}$. The case $\ba\ne\ba^{*_\g}$ now
follows as in the proof of Theorem~3.

At last we are ready to prove Theorem~2. Recall the Janssen
representation \eqref{12}, which is valid for the frame operator
$S_{\phi,\phi}$ under consideration by Proposition~2. In terms of
the notation of Lemma~4, it says that
\begin{equation*}
S_{\phi,\phi}=\pi(\bb), \quad\mbox{where} \ \ b_{jk}=\langle \phi,
M_{j/\b}T_{k/\a}\phi\rangle.
\end{equation*}
Since $\phi\in M^1_v$, it is easy to verify that $\bb\in
l^1_{\tilde v}$, and since $\CG(\phi,\a,\b)$ is a frame for $L^2$,
$S_{\phi,\phi}=\pi(\bb)$ is invertible on $L^2$. Hence $\bb$ is
invertible in $l^1_{\tilde v}$ by Lemma~8, so
$S_{\phi,\phi}^\inv=\pi(\bb^\inv)$ is bounded on $M^1_v$ by
Lemma~5, and we are done.

\end{document}